
\hsize=6.5 true in
\vsize=8.7 true in
\baselineskip=15pt
\vglue 1.5 cm


\def\qed{$\rlap{$\sqcap$}\sqcup$}

\def\a{\bigskip \par \noindent}
\def\b{\par \noindent}
\font\maiu=cmcsc10
\font\small=cmr9 at 9truept

\font\tengothic=eufm10
\font\sevengothic=eufm7
\newfam\gothicfam
       \textfont\gothicfam=\tengothic
       \scriptfont\gothicfam=\sevengothic
\def\goth#1{{\fam\gothicfam #1}}


    \font\tenmsb=msbm10              \font\sevenmsb=msbm7
\newfam\msbfam
       \textfont\msbfam=\tenmsb
       \scriptfont\msbfam=\sevenmsb
\def\Bbb#1{{\fam\msbfam #1}}


\def\PP#1{{\Bbb P}^{#1}}

\def\PPP {{\Bbb P}}

\def\ref#1{[{ #1}]}
\vskip 1cm 
\centerline{\bf On the minimal free resolution for fat point schemes}
\centerline{\bf of multiplicity at most $3$ in $\Bbb {P}^2$.} 
\bigskip 
\centerline{ Edoardo Ballico, Monica Id\`a}
\footnote{}{\rm The authors were partially supported by MUR and GNSAGA of INdAM (Italy).} 

\bigskip 
\noindent 
{\small {\maiu Abstract.} Let $Z$ be a fat point scheme in $\Bbb {P}^2$ supported on general points. Here
we prove that if the multiplicities are at most 3 and the length of $Z$ is sufficiently high then the number of generators of the homogeneous ideal $I_Z$ in each degree is as small as numerically possible. Since it is known that $Z$ has maximal Hilbert function, this implies that $Z$ has the expected minimal free resolution.}

\smallskip \noindent {\maiu msc: 14N05} 
\b {\small {\maiu keywords:} minimal free resolution; homogeneous ideal;
zero-dimensional scheme; fat double point; fat triple point; projective plane.}
\bigskip 
\noindent{\bf 1. Introduction.} 
\medskip 

What is the minimal free resolution of a `` general ''
zero-dimensional scheme $Z\subset \Bbb {P}^2$? In this paper ``
general '' means that $\sharp (Z_{red})$ is fixed and $Z_{red}$ is general in $\Bbb {P}^2$. We are interested in the minimal free resolution of general fat point schemes of
$\Bbb {P}^2$. A fat point $mP\subset  \Bbb {P}^2$ is the zero-dimensional
subscheme of $\Bbb {P}^2$ with support in the point $P$ and $({\cal I}_P)^{m}$ as its ideal
sheaf. A general fat point scheme $Z:= m_1P_1+\dots +m_rP_r$ of $\Bbb {P}^2$, $m_1\geq \ldots m_r\geq 0$, is a general zero-dimensional scheme such that for each $P_i\in Z_{red}$ the connected component of $Z$ with support in $P_i$ is the fat point $m_iP_i$. If $m=1$ (resp. $m=2$, resp. $m=3$) we will say that $mP$ is a simple (resp. double, resp. triple) point.
We recall that ${\rm length}(mP) =
m(m+1)/2$ for all $m>0$. 
\par The Hilbert function and the minimal free resolution of plane
fat point schemes  have been studied quite a lot in the last years, assuming that the number of points is low, or that the multiplicities are low, or giving some other kind of restriction on the involved integers. For example, the Hilbert function is known in the equimultiple case for any $r$ if $m=m_1=\ldots =m_r\leq 20$ (\ref {Hi}, \ref {C-C-M-O}) and in many other cases (\hbox{\ref {H-H-F}}, \ref {H-R}, \ref {E2}, \ref {R}); it is also known if $r\leq 9$ (\ref {N}, \ref {Ha2}) and if $m_i\leq 7$ (\ref {M}, \ref {Y}). The graded Betti numbers for a fat point scheme $Z=m_1P_1+\dots +m_rP_r$ are known if $r\leq 8$ (\ref {Ca}, \ref {F}, \ref {Ha3}, \ref {F-H-H}) and in some other cases (\ref {H-H-F}, \ref {H-R}). For the equimultiple case, there is a general conjecture (\ref {Ha1}), proved for $m\leq 3$; i.e., it is known that the homogeneous ideal $I(Z)$ of $Z$ is minimally generated for $m=1$ (\ref {G-M}), $m=2$ (\ref {I})
and $m=3$ (\ref {G-I}). 
\par If $Z$ is a fat point scheme of multiplicity at most $3$ in $\Bbb {P}^2$, $Z$ has maximal Hilbert function by \ref M, i.e. $h^1(\Bbb {P}^2, {\cal I}_Z(k)) \cdot h^0(\Bbb {P}^2,{\cal I}_Z(k)) = 0$ for all $k \ge 0$, provided that $k\geq m_1+m_2+m_3$. In this paper we show that, if the length of $Z$ is sufficiently high, the multiplication maps $\mu _k(Z): H^0({\cal I}_{Z}(k))\otimes H^0({\cal O}_{\PP 2}(1)) \rightarrow H^0({\cal I}_{Z}(k+1))$ are of maximal rank for any $k$. The following result hence follows:

\medskip {\bf Theorem 1.1.} {\it Fix non-negative integers $a, b, c$ such that $a+3b+6c \ge 79$, and let $Z \subset \Bbb {P}^2$ be a general union of $\,a$ simple points, $\,b$ double points and $\,c$ triple points in $\Bbb {P}^2$. Let
$v$ be the minimal integer such that $a+3b+6c \le (v+2)(v+1)/2$.
Then $Z$ has the expected
minimal free resolution, i.e. the homogeneous ideal
of $Z$ is minimally generated by $(v+2)(v+1)/2 -a-3b-6c\,$ forms of
degree $v$ and $\max \{0,2a+6b+12c-v^2-2v\}$ forms
of degree $v+1$.}

\medskip  In the case in which there are only double points, we can handle a
few low integers $v$ and prove the following complete result:

\medskip {\bf Theorem 1.2.} {\it Fix integers $a, b$ such that $a \ge 0$, $b \ge 0$, and let $Z \subset \Bbb {P}^2$ be a general union of $\, a$ simple points
and $\, b$ double points of $\, \Bbb {P}^2$. Let
$v$ be the minimal integer such that $a+3b \le (v+2)(v+1)/2$.
Then, for any $(a,b)\neq (0,2),(0,5),(1,1),(1,2)\,$,  $Z$ has the expected
minimal free resolution, i.e. the homogeneous ideal
of $Z$ is minimally generated by $(v+2)(v+1)/2 -a-3b\,$ forms of degree
$v$ and $\max \{0,2a+6b-v^2-2v\}$ forms
of degree $v+1$.}

\par We raise the following conjecture.

\medskip {\bf Conjecture 1.3.}
Fix integers $m>0$ and $n\ge 2$. Then there is an integer $\alpha _{n,m}$
such that for all integers $r> \alpha _{n,m}$
the minimal free resolution
of a general union $Z \subset \Bbb {P}^n$ of $r$ multiple points
with multiplicity at most $m$ is the expected one.

\medskip For related conjectures and discussions, see also \ref {Ha1} and \ref {H-H-F}. Also recall the Minimal Resolution Conjecture, which is about the generic union of $r$ simple points in $\PP n$, and has been proved for $n=2$ (as said above), for $n=3$ (\ref {B}, \ref {B-G}), in particular cases for $n\geq 4$, and for any $n$ and $r>>0$ (\ref {H-S}).

\medskip Several proofs (all of them heavily using the Horace method) of Theorem 1.1 and Theorem 1.2 might be given. We will see in section 2 that Theorem 1.2 is an obvious consequence of \ref {I} and of the deformation of a double point to three simple points (see Remark 2.2). The proof of Theorem 1.1  is obtained as well adapting the proof of Theorem 1.1 in \ref {G-I}, but it needs a bigger effort. 
\par Notice that there is no serious obstacle in proving Theorem 1.1 when length $Z$ is lower and $c\neq 0$, it is essentially a question of patience; for example, in order to go down to $a+3b+6c\geq 56$ we only need to consider a few more cases, thanks to Lemma 3.1.5 and Lemma 3.1.6; if the fat point scheme is supported on $\leq 8$ points, \ref{F-H-H} give the Betti numbers; if $a,b,c$ are such that $a+3b\equiv 0$  (mod 6), $3b\leq a$ and $\lfloor {a+3b\over 6}\rfloor\neq2,3,5$, then \ref {G-I} and Remark 2.2 gives the result through semicontinuity. There are low length cases where the minimal free resolution is not the expected one, for example the union of 2 or 3 or 5 triple points. 

\par In the following we use the Horace method, introduced by A.Hirschowitz to prove this interpolation type problems, and the differential Horace method. For the first one, we send the reader to \ref {Hi}. The second one has been introduced in \ref {A-H} for invertible sheaves, and then extended to vector bundles (see \hbox{\ref {G-I}} Proposition 2.6).  
\par We work over an algebraically closed field $ \Bbb K$ such that either $\hbox{char}(\Bbb {K}) = 0$ or
$\hbox{char}(\Bbb {K})>3$. In \ref {I} and \ref {G-I} there is the
$\hbox{char}(\Bbb {K}) = 0$ assumption, but in fact
this assumption is not necessary in \ref {I} , while in \hbox{\ref {G-I}}  it is used only
in the proof of Lemma 2.9 where a map $t_i \mapsto t_i^{r_i}$ between
formal power series rings is considered. To get the injectivity of the differential
of this map at $(0,\dots ,0)$ it is sufficient to assume $\hbox{char}(\Bbb {K}) > r_i$ for all
$i$. Since in our set up (fat point with multiplicity at
most $3$) we have $r_i \le 3$ for all $i$, we hence assume $\hbox{char}(\Bbb {K}) = 0$ or
$\hbox{char}(\Bbb {K})>3$.

{\maiu Acknowledgements.} The authors wish to thank the referee for her/his accurate report and useful remarks. 

\a {\bf 2. Preliminaries and proof of Theorem 1.2.}

\a {\bf 2.1.} Let $X$ be a 0-dimensional scheme in $\PP 2$, and let length$X=l$; we say that $X$ has maximal Hilbert function in degree $k$ if $h^0({\cal I}_X(k))=max\{0, {k+2\choose 2}-l\}$; $X$ has maximal Hilbert function if this is true for any $k$.
\par Now assume that $X$ has maximal Hilbert function. Then $X$ has the expected minimal free resolution if the natural multiplication maps
$$ \mu _k: H^0({\cal I}_X(k))\otimes H^0({\cal O}_{\PP 2}(1)) \rightarrow H^0({\cal I}_X(k+1)) $$ 
have maximal rank for each $k$.
Set $\Omega:=\Omega^1_{\PP 2}$. Taking the cohomology sequence of the Euler sequence in $\PP 2$ tensored by the ideal sheaf ${\cal I}_X(k+1)$, we see that
ker $\mu _k= H^0(\Omega (k+1)\otimes {\cal I}_X)$. 
\par Set $v=v(X)=min \{k\geq 1\ \vert\ {k+2\choose 2}-l\geq 0\}$; then, $\mu_k$ is trivially injective for $k< v$ because $h^0({\cal I}_X(k))=0$, and $\mu_k$ is surjective for $k>v$ by the Castelnuovo-Mumford Lemma because $h^1({\cal I}_X(v))=0$.
Hence $X$ is minimally generated if and only if $\mu_v$ is of maximal rank.

\medskip Now let $w=w(X)=min \{k\geq 1\ \vert\ k(k+2)-2l\geq 0\}$ ($w$ is the smallest integer for which the restriction map $\rho_k: H^0(\Omega(k+1))\to H^0(\Omega(k+1)\vert _X)$ can be surjective). Then $2l\leq w(w+2)$ gives $l< {w+2\choose 2}$, and by assumption $X$ has maximal Hilbert function, hence $h^1({\cal I}_X(w))=0$ and $h^0({\cal I}_X(w))={w+2\choose 2}-l>0$. Now if $h^1({\cal I}_X(k))=0$ then $h^1({\cal I}_X(k+1))=0$ and we get $3h^0({\cal I}_X(k))- h^0({\cal I}_X(k+1))=3({k+2\choose 2}-l)-({k+3\choose 2}-l)=k(k+2)-2l$. Hence $w$ is the smallest integer $k$ for which the map $\mu_k$ can be surjective, without being the 0-map. From what is said above we have  $v\leq w\leq v+1$, so that $X$ is minimally generated if and only if $\mu_{w-1}$ is injective and $\mu_{w}$ is surjective, and this happens if and only if $h^0(\Omega (w)\otimes {\cal I}_X)=0$ and $h^0(\Omega (w+1)\otimes {\cal I}_X)=w(w+2)-2l$. 

\medskip If we assume only that the Hilbert function of $X$ is maximal in degree $w$, without any assumptions on the other degrees, then the same considerations as above show that $\mu_{w-1}$ is injective if and only if $h^0(\Omega (w)\otimes {\cal I}_X)=0$, and $\mu_{w}$ is surjective if and only if $h^0(\Omega (w+1)\otimes {\cal I}_X)=w(w+2)-2l$.

\a What we do in practice is to look for suitable schemes for which the map $\rho_k$ is bijective, and from these deduce the injectivity or surjectivity for the schemes we are interested in. For arithmetical reasons (think of $k$ odd) it is better to work in the projectivized bundle $\PPP (\Omega)$ with the canonical projection $\pi: \PPP (\Omega) \to \PP 2$. We set ${\cal E}_{k}:={\cal O}_{\PPP (\Omega)}(1)\otimes \pi^* {\cal O}_{\PP 2}(k)$. One has (e.g., see \ref {I2} Lemma 2.1 ) $$H^0(\Omega
(k+1)\otimes {\cal I}_X)\cong H^0({\cal E}_{k+1}\otimes {\cal I}_{\pi^{-1}X})\,,\quad  \quad H^0(\Omega
(k+1)\vert_X)\cong H^0({\cal E}_{k+1}\vert_{\pi^{-1}X}).$$ 
\medskip If $X\subset \PPP (\Omega)$ is a 0-dimensional scheme such that $length X=H^0({\cal E}_{k+1})$, and $H^0({\cal E}_{k+1}\otimes {\cal I}_{X})=0$, we say that $X$ is $k$-settled.

\a {\bf Remark 2.2.} It is immediate to see that a $2$-fat point of $\PP 2$ is the flat limit of a family whose general fiber is the general union of 3 simple points. 
\par A $3$-fat point of $\PP 2$ is the flat limit of a family
whose general fiber is the general union of one $2$-fat point and $3$ simple points (see \ref {E}, Proposition 4).
\par Hence, a $3$-fat point of $\PP 2$ is the flat limit of a family whose general fiber is the general union of 6 simple points.

\a {\bf Notations 2.3.} We denote by $R_p$, with $p=0,1,2,3,5,8,11$, a certain
0-dimensional scheme of length $p$ in $\PPP (\Omega)$ which we now define.  Let $U$ be an open subset in $\PP 2$ and $\Omega\vert _U \cong E_1\oplus E_2$ a local trivialization for $\Omega$; then, if $p\neq 0,1$, $R_{p}= \eta_1\cup
\eta_2$, where $\eta_1$, $\eta_2$ have support on two distinct points $A_1,A_2$
in the same fiber $\pi ^{-1}(P)$ where
$A_i=\PPP (E_i)\cap \pi ^{-1}(P)$ and
$\eta _i \subset \PPP (E_i)$, $i=1,2$, so that ${\rm length} (\eta_i\cap \pi
^{-1}(P))=1$. If we consider affine coordinates
$\{ x,y,z\}$ in an affine chart of $\PPP (\Omega)$ containing $R_{p}$, we may
suppose $\pi ^{-1}(P)=\{x=y=0\}$,
$A_1=(0,0,0)$ and $A_2=(0,0,1)$; then $R_{p}$ is defined as follows:
\b $R_0 = \emptyset$; $R_1=\{A_1\}$ is just a point in $\PPP (\Omega)$; $R_2=
\{A_1,A_2\}$ ;
\b $R_3$ is made of the point $\eta_1=A_1$ and a length 2
structure $\eta_2$ on $A_2$, given by an ideal of type $(x,y^2,z-1)$;
\b $R_5$ is made of a length 2 structure $\eta_1$ on $A_1$, given by an ideal of
type $(x,y^2,z)$ and the first infinitesimal neighbourhood $\eta_2$ on $A_2$,
given by an ideal of type $(x^2,xy,y^2,z-1)$;
\b $R_8$ is made of two 4-ple structures $\eta_1$, $\eta_2$ of the same type,
given by ideals of type
$(x^3,xy,y^2,z)$, $(x^3,xy,y^2,z-1)$;
\b $R_{11}$ is such that $\eta_1$ is a 5-ple structure on $A_1$ given
by an ideal of type $(x^3,x^2y,y^2,z)$, and $\eta_2$
is given by an ideal of type $(x^3,x^2y,y^2x,y^3,z-1)$.
\a So for $k$ even $R_p$ is two copies of a nilpotent $\eta\subset \PP 2$ with
$\eta\cong \eta_i$, while for $k$ odd $R_p$ is the same thing but with a
nilpotent transversal to one of the components of this scheme added. Since we are interested in the
schemes $R_p$  only for the vanishing of global sections of ${\cal O}_{\PPP
(\Omega)}(1)\otimes \pi^*{\cal O}_{\PP 2}(t)$, we can consider (see Lemma 2.2
in \ref {G-I}) that $R_p$ is the pull back of $\eta \subset \PP 2$ for $k$ even, and
for $k$ odd the pull back of $\eta $ with a nilpotent transversal to this scheme
added.
For the same reason, if $B \subset U$ is any 
0-dimensional subscheme of $\PP 2$ and we set: 
\b \centerline{$ B':= \pi^{-1}(B)\cap \PPP (E_1), \qquad  B^{\prime \prime}:= 
\pi^{-1}(B)\cap \PPP (E_2), \qquad 
\widehat{B}:= B'\cup B^{\prime \prime} $,}
\b as long as we are concerned only with the vanishing of the global section of ${\cal E}_k$ along $\pi^{-1}(B)$, we can substitute the last one with $\widehat{B}$. 

\a {\bf Notations 2.4.} With $Y(a,b)$ we denote in the following the generic union of $a$
points and $b$ double points in $\PP 2$; we also set $\tilde Y(a,b):=\pi^{-1}(Y(a,b))$.
\b For any $k\geq 0$, let $q=q(k)$, $r=r(k)$ be positive integers such that
$k(k+2)=6q(k)+r(k)$, with $0\leq r \leq 5$; the possible values for $r$ are
$0,2,3,5$ (see \ref {I} Lemma 1.6).
\b In the following $Z(s,d,p)$ will denote the generic union in $\PPP (\Omega)$
of $\tilde Y(s,d)$ with $R_p$, where $2s+6d+p=k(k+2)$ and $0\leq p \leq r$. Notice that $p\equiv r$ (mod 2), hence $p\neq 4$.
\b We set  $\Delta _k=\{(s,d,p)\in \Bbb N ^3 \vert \; 2s+6d+p=k(k+2), \;
p=0,1,2,3,5 \}$.
\b In \ref {I} the assertion:``$Z(0,q(k),r(k))$ is k-settled", denoted by ``{\it A(k)}", is proved for any $k\neq 2,3$.

\a {\bf Lemma 2.5.} {\it If A(k) in \ref {I} is true and if $(s,d,p)\in \Delta _k$,
then $H^0({\cal E}_{k+1}\otimes
{\cal I}_{Z(s,d,p)})=0$.}
\a {\bf Proof.} We write $q=q(k)$, $r=r(k)$; by assumption $2s+6d+p=6q+r$ with $0\leq p\leq r\leq 5$. Writing $s=3l+j$, $0\leq j \leq 2$, we find $4+r\geq 2j+p \equiv r$ (mod 6), hence $2j+p=r$.
\b Now observe that a double point is specialization of 3 points in the plane;
$R_2$ is the pull back of a point of $\PP 2$; $R_3$ is specialization of the
union of the pull back of a point of $\PP 2$ with $R_1$, which is a point in
$\PPP (\Omega)$, and if $2j+p=5$, then $R_5$ is specialization of the general
union of the pull back of $j$ points of $\PP 2$ with $R_p$ ($j=1,2$).
\b Hence, the scheme $Z(s,d,p)$ specializes to $Z(0,q,r)$ and we conclude by semicontinuity since by assumption $H^0({\cal E}_{k+1}\otimes {\cal I}_{Z(0,q,r)})=0$. \qed

\a {\bf Lemma 2.6.} {\it Let $a,b,k$ be nonnegative integers such that $(k-1)(k+1) <2a+6b\leq
k(k+2).$
 Then,
\b if A(k-1) is true, $\mu _{k-1}$ is injective for $Y(a,b)$;
\b if A(k) is true, and if $Y(a,b)$ has maximal Hilbert function in degree $k$, then $\mu _{k}$
is surjective for $Y(a,b)$.}

\a {\bf Proof.} By assumption, $w(Y(a,b))=k$. Then, in order to prove the first and respectively the second statement, it is enough to prove (see preliminaries 2.1) $h^0({\cal E}_{k}\otimes {\cal I}_{\tilde Y(a,b)})=0$ and respectively
$h^0({\cal E}_{k+1}\otimes {\cal I}_{\tilde Y(a,b)})=k(k+2)-2(a+3b)$.

\par Assume A(k-1) holds.
We now show that there exists $(s,d,p)\in \Delta _{k-1}$ such that $Z(s,d,p)$ is
contained in $\tilde Y(a,b)$. We write $q=q(k-1)$, $r=r(k-1)$; we are looking for
$s,d,p$ such that $2s+6d+p=(k-1)(k+1)=6q+r<2a+6b$, with $0\leq r \leq 5$. Set
$r=2l+\epsilon$, $\epsilon =0,1$. If $q<b$, we set $d=q$, $s=0$, $p=r$; since
$R_r$ is contained in the pull back of a double point, we have $Z(0,q,r)\subset
\tilde Y(a,b)$. If $q\geq b$, we set $d=b$, $s=3(q-b)+l$, $p=\epsilon$. Since
$2(3(q-b)+l)+\epsilon < 2a$, we have $s+1\leq a$ and
moreover if $\epsilon =1\,$ $R_\epsilon$ is
contained in the pull back of a simple point, hence
$Z(3(q-b)+l,b,\epsilon)\subset \tilde Y(a,b)$ both for $\epsilon =0$ and for $\epsilon =1$.

\smallskip By the previous lemma  $H^0({\cal E}_{k}\otimes {\cal
I}_{Z(s,d,p)})=0$; we conclude taking cohomology of the exact sequence:
$0\rightarrow {\cal E}_{k}\otimes {\cal I}_{\tilde Y(a,b)}\rightarrow {\cal
E}_{k}\otimes {\cal I}_{Z(s,d,p)}\rightarrow...$.

\a Now assume A(k) holds.
We now show that there exists $(s,d,p)\in \Delta _{k}$ such that $Z(s,d,p)$
contains $\tilde Y(a,b)$. We write $q=q(k)$, $r=r(k)$; we are looking for $s,d,p$ such
that $2a+6b\leq 2s+6d+p=k(k+2)=6q+r$, with $0\leq r \leq 5$. Set
$r=2l+\epsilon$, $\epsilon =0,1$. We have $6(b-q)\leq r-2a\leq 5$, hence $b\leq
q$, so we set $d=b$, $s=3(q-b)+l$, $p=\epsilon$; since $2a\leq
2(3(q-b)+l)+\epsilon $, we have $ a\leq s$, so $Z(3(q-b)+l,b,\epsilon)\supset
\tilde Y(a,b)$.

\a By the previous lemma  $H^0({\cal E}_{k+1}\otimes {\cal I}_{Z(s,d,p)})=0$;
the cohomology of the exact sequence:

$$0\rightarrow {\cal E}_{k+1}\otimes {\cal I}_{Z(s,d,p)}\rightarrow {\cal
E}_{k+1}\otimes {\cal I}_{\tilde Y(a,b)}\rightarrow {\cal E}_{k+1}\otimes {\cal
I}_{\tilde Y(a,b),Z(s,d,p)}\rightarrow 0$$

\b gives $h^0({\cal E}_{k+1}\otimes {\cal I}_{\tilde Y(a,b)})\leq h^0({\cal E}_{k+1}\otimes
{\cal I}_{\tilde Y(a,b),Z(s,d,p)})=2(s-a)+6(d-b)+p=k(k+2)-(2a+6b)$. On the other hand $h^0({\cal
E}_{k+1}\otimes {\cal I}_{\tilde Y(a,b)})\geq h^0({\cal E}_{k+1})-h^0({\cal
E}_{k+1}\vert _{\tilde Y(a,b)})=k(k+2)-(2a+6b)$, so we have equality. 
 \qed

\a {\bf Proof of Theorem 1.2.} First let us check that $Y(a,b)$ has maximal Hilbert function (mHf for short) for $(a,b)\neq (0,2),
(0,5)$. It is well known that $Y(a,b)$ has mHf for any $a$ if $b=0$, or for any $b\neq 2,5$ if $a=0$ (\ref{Hi}). Now let $\goth L \neq 0$ be a linear system in $\PP 2$; if $P$ is a point outside of the base locus of $\goth L$, for example a generic point of the plane, and $\goth L (P)$ is the linear system obtained by $\goth L$ imposing the passage through $P$, then $\rm dim \goth L (P)=  \rm dim \goth L-1$. Hence if we add (generic) simple points to a scheme $Y(a,b)$ with mHf we get a scheme which again has mHf. We conclude that all schemes $Y(a,b)$ with $b\neq 2,5$ have mHf. Moreover, if $a\geq 3$, and $b=2,5$, the scheme $Y(a,b)$ specializes to $Y(a-3,b+1)$ which has 
mHf, hence by semicontinuity $Y(a,b)$ has mHf too. It is immediate to check by hand that also in the remaining cases, i.e. $(a,b)=(1,2),(2,2),(1,5),(2,5)$ the  Hilbert function is maximal.

\medskip Now we study the maps $\mu_k$. Recall that assertion $A(k)$ is proved in \ref{I} for $k=1$ and $k\geq 4$. If $l=l(Y(a,b))=a+3b>12$, and $w$ is the integer such that $(w-1)(w+1) <2a+6b\leq w(w+2)$, then $w=w(Y(a,b))\geq 5$ and Lemma 2.6 assures that $Y(a,b)$ is minimally generated, provided that $(a,b)\neq (0,5)$ (see also 2.1). If $w=4$ (i.e. $15<2l\leq 24)$ and $v=v(Y(a,b))=4$ (i.e. $10<l\leq 15)$, then (see 2.1) it is enough to prove that $\mu_4$ is surjective, and this is true by Lemma 2.6; hence $Y(a,b)$ is minimally generated also for $11\leq l\leq 12$. 

\par We now assume $l\leq 10$. We know that $Y(a,b)$ is minimally generated for any $a$ if $b=0$ (\ref{G-M}), or for any $b\neq 2$ if $a=0$ (\ref{I}); the remaining cases are $\{(a,1)\}_{a=1,\dots,7}$, $\{(a,2)\}_{a=1,\dots,4}$ and $(1,3)$.
For $(7,1),(4,2),(1,3)$ we have $h^0({\cal I}_{Y(a,b)}(3))=h^1({\cal I}_{Y(a,b)}(3))=0$ so by the Castelnuovo-Mumford Lemma the scheme is minimally generated. The scheme $Y(6,1)$ specializes to $Y(3,2)$ which specializes to $Y(0,3)$, and the last one is minimally generated, hence by semicontinuity the other two are minimally generated (notice that all schemes here have mHf, and this is why semicontinuity for $H^0(\Omega (k+1)\otimes {\cal I})$ is useful). 
\par The few cases left can be recovered from \ref{Ca}; anyway, we give their explicit description in what follows. For $(3,1),(2,1)$ we have $v=2$ and $w=3$, so it is enough to prove that $\mu_2$ is injective, and this is true because $h^0({\cal I}_{Y(3,1)}(2))=0$ and $h^0({\cal I}_{Y(2,1)}(2))=1$.
For $(2,2)$ we have $v=3$ and $w=4$, so it is enough to prove that $H^0(\Omega (4)\otimes {\cal I}_{Y(2,2)})=0$. Let $C$ be a conic through the four points, and $L$ the line through the two double points; since $\Omega \vert_C\cong {\cal O}_{\PP 1}(-3)^{\oplus 2}$, we have $H^0(\Omega (4)\vert_C\otimes {\cal I}_{Y(2,2)\cap C,C})=0$, hence $h^0(\Omega (4)\otimes {\cal I}_{Y(2,2)})=h^0(\Omega (2)\otimes {\cal I}_{Y(2,0)})$ and the last one is zero because $H^0(\Omega (1))=0$ and $H^0(\Omega (2)\vert_L\otimes {\cal I}_{Y(2,0),L})=0$. The scheme $Y(5,1)$ specializes to $Y(2,2)$, hence also $Y(5,1)$ is minimally generated.
\par For $(4,1)$ we have $v=3$ and $w=3$, so we want to prove that $\mu_3$ is surjective, or equivalently (see 2.1) that $ h^0({\cal E}_{4}\otimes {\cal I}_{\tilde Y(4,1)})=h^0(\Omega (4)\otimes {\cal I}_{Y(4,1)})=15-2l=1$; this is in turn equivalent to $ H^0({\cal E}_{4}\otimes {\cal I}_{Z(4,1,1)})=0$, where $Z(4,1,1)=\tilde Y(4,1)\cup R_1$. Now we use the Horace method in $\Bbb P (\Omega)$, as we do, for example, in the proof of the forthcoming Lemma 3.1.4.
Let $C$ be the conic through the five points in the support of $Y(4,1)$; then, $H^0({\cal E}_{4}\vert_{\pi^{-1}C}\otimes {\cal I}_{Z(4,1,1)\cap \pi^{-1}C,\pi^{-1}C)})=0$, and $H^0({\cal E}_{2}\otimes {\cal I}_{Z(1,0,1)})=0$ (the last one is assertion $A(1)$), so we conclude that $ H^0({\cal E}_{4}\otimes {\cal I}_{Z(4,1,1)})=0$.

\par $Y(1,1)$ is not minimally generated. In fact, let $L$ be the line through the two points. Here $v=2$ and $w=2$, but $\mu_2$ cannot be surjective since $L$ is in the base locus of $H^0({\cal I}_{Y(1,1)}(2))$, so there must be a generator of degree 3.
$Y(1,2)$ is not minimally generated. Here $v=3$ and $w=3$, but $\mu_3$ cannot be surjective since the line through the two double points is in the base locus of $H^0({\cal I}_{Y(1,1)}(3))$, so there must be a generator of degree 4. \qed

\a {\bf 3. \bf Proof of Theorem 1.1.}

\a { \maiu 3.1 Reduction to a statement with no simple point.}

\a {\bf 3.1.1. Recall of techniques and notations from \ref {G-I}.} In the following we use, beyond the Horace method, also  the differential Horace method for vector bundles, for which we refer to \hbox{\ref {G-I}} Section 2 and in particular Proposition 2.6.; moreover, we'll use notations 3.1 and the ones estabilished 
at the beginning of the proof of 3.3 in \hbox{\ref {G-I}}, so we recall them briefly here.  

\par For the definition of vertically graded subscheme with base a fixed smooth divisor see \ref {A-H} and \ref {G-I} 2.3. 
We introduce now some notations that will allow us to express 
ourselves as if we were working in 
$\PP 2$, while our environment is actually $\PPP (\Omega)$. 
Let $B$ be a 0-dimensional scheme of $\PP 2$ with support at a point 
$P$, vertically graded  with base a smooth conic $C$ with local 
equation $y=0$, and let $x,y$ be local coordinates at $P$. Notice that $B'$ and $B''$ (see Notations 2.3) are vertically graded with base $H=\pi ^{-1}(C)$. Consider the integers $a_j$ where, if ${\cal I}_{Tr^j_C(B)} := ({\cal I}_{B}:{\cal 
I}_C^j)\otimes{\cal O}_C$, then $Tr^j_C(B)=(x^{a_j},y)$.
We will denote $\widehat B=B'\cup B''$ by $\pmatrix {a_s\cr \vdots \cr a_0}$. So for example $\widehat B$ is denoted by $(h)$  if $I_B=(x^h,y)$; by $\pmatrix {1 \cr 2}$ if $I_B=(x,y)^2$; by $\pmatrix {1 \cr 3}$ if $I_B=(x^3,xy,y^2)$; by $\pmatrix {2 \cr 3}$ if $I_B=(x^3,x^2y,y^2)$; by $\pmatrix {2 \cr 2}$ if $I_B=(x^2,y^2)$; by $\pmatrix {1 \cr 2 \cr 3}$ if $I_B=(x,y)^3$.

\par If $\widehat B$ is a $\pmatrix {1 \cr 3}$, we say that $\widehat B$ is ``a 
$\pmatrix {1 \cr 3}$ scheme over $C$" and write $\pmatrix {1 \cr 3}_{over C}$.  We write, 
for example, ``$h\pmatrix {1 \cr 2}$ schemes general over $\PP 2$ " to mean 
a union of $h$ schemes of type $\pmatrix {1 \cr 2}$, whose projection in 
$\PP 2$ is general. 
Moreover, if  $h, l\in { \Bbb N}$, we will use, for example, the 
notation  ``$h\pmatrix {1 \cr 2}+ l\pmatrix {1 \cr 3}$" 
to denote the union of $h$ schemes of type $\pmatrix {1 \cr 2}$ and  $l$ 
schemes of type $\pmatrix {1 \cr 3}$.
\par If $a_i$ and $b_i$, $i=1,...,m$ are positive integers, with an abuse of notation we write
$(\sum _{i=1}^m a_ib_i)=\sum _{i=1}^m a_i(b_i)$, since for our vanishing problem only the length of the scheme over $C$ matters. 
\a When we apply \ref {G-I} Proposition 2.6, i.e. when we do a differential Horace, or HD, step, and we say for example that we ``add over $C$"  $\; \pmatrix {1 \cr [3]}+\pmatrix {1 \cr [2]\cr 3}$, this means that we are using the ``ground slide" of the vertical scheme $\pmatrix {1 \cr 3}$ and the ``first floor slide" of the triple point, so that the HD trace $\pmatrix {5}_C$ is given by the numbers in square brackets, while the HD residue $(\pmatrix {1 }+\pmatrix {1 \cr 3})_C$ is obtained by the eliminating the part in the square brackets.

\a {\bf Notations 3.1.2.} With $Y(a,b,c)$ we denote in the following the generic union of $a$
points, $b$ double points and $c$ triple points in $\PP 2$; we also set $\tilde
Y(a,b,c):=\pi^{-1}(Y(a,b,c))$.
\b For any $k\geq 0$, let $u=u(k)$, $\rho =\rho (k)$ be the positive integers
such that $k(k+2)=12u+\rho$, with $0\leq \rho \leq 11$. 
If we write  $k$ modulo 6, we get (see \ref {G-I} 1.2; but notice that there $u(k)$, resp.$\rho (k)$, are denoted by $q(k)$, resp.$r(k)$):

\item{} for k=6l $\; \;\;$ \qquad   $u(k)=3l^2+l$  $\; \; \;\;\;\,\,$ \qquad $\rho (k)=0$
\item{} for k=6l+1 \qquad  $u(k)=3l^2+2l$  $\; \; \;\;\;\,$\qquad $\rho (k)=3$  
\item{} for k=6l+2 \qquad  $u(k)=3l^2+3l$  $\; \; \;\;\;\,$\qquad $\rho (k)=8$  
\item{} for k=6l+3 \qquad  $u(k)=3l^2+4l+1$ \qquad $\rho (k)=3$  
\item{} for k=6l+4 \qquad  $u(k)=3l^2+5l+2$ \qquad $\rho (k)=0$  
\item{} for k=6l+5 \qquad  $u(k)=3l^2+6l+2$ \qquad $\rho (k)=11$.

\bigskip For a fixed $k$, let $Z(s,d,t,p)$ denote the generic union in $\PPP (\Omega)$ of $\tilde
Y(s,d,t)$ with $R_p$, where $2s+6d+12t+p=k(k+2)$ and $0\leq p\leq \rho$, $p=
0,1,2,3,5,8,11$.
\b We'll set $\Lambda_k=\{ (s,d,t,p)\in  \Bbb N ^4 \vert \;
2s+6d+12t+p=k(k+2),\; 0\leq p\leq \rho ,\, p= 0,1,2,3,5,8,11\}$.
\b In the following, $``H(s,d,t,p,k)"$ denotes the statement:
\medskip \centerline {``If $(s,d,t,p)\in \Lambda _k$, then $H^0({\cal E}_{k+1}\otimes
{\cal I}_{Z(s,d,t,p)})=0.$"}
\a We want to prove that
$H(s,d,t,p,k)$ holds for $k\geq 12$, and this will be done through some lemmas.
Notice
that if $t=0$ the statement $H(s,d,t,p,k)$ reduces to Lemma 2.5, so we can assume
$t\geq 1$.

\a {\bf Remark 3.1.3.}  Recall that $B(k)$ in \ref {G-I} is nothing else but $H(0,0,u(k),\rho (k),k)$,
and $B(k)$ is true for $k\geq 10$ (\ref {G-I} Proposition 3.9). To be punctual about that, notice that for the remainder scheme $T_8$ defined in \ref {G-I} 1.3, which is the analogous of our scheme $R_8$, one has the choice between two 4-ple structures $\eta_1$, $\eta_2$ of the same type, given in local coordinates $x,y$ by an ideal of type $(x^3,xy,y^2)$ or by an ideal of type $(x^2,y^2)$; in the languages of vertical schemes, $T_8$ is $\pmatrix { 1\cr 3}$ or $\pmatrix { 2\cr 2}$. Anyway, the unique point of the proof of Theorem 1.1 in \ref {G-I} where one chooses to use $T_8$ as a $\pmatrix { 2\cr 2}$ is the next to last step in the proof of Proposition 4.1, where we specialize it on (the pull back of ) a smooth conic $C$. It is possible to choose $T_8=\pmatrix { 1\cr 3}$ also here, specializing it as a $\pmatrix { 1\cr 1\cr [2] }$ (in other words, we consider it as vertical scheme with respect to the $y$-axis instead that to the $x$-axis); in the last step, the residue $\pmatrix { 1\cr 1}$ can be specialized as a $(2)$ on $C$. Hence we can always assume that the remainder scheme $T_8$ in \ref {G-I} is our scheme $R_8$.

\bigskip The following lemma  allows to construct a lot of well generated schemes
containing  2-fat points; the proof goes exactly as the  proof of Lemma 2.1 in
\ref {I}, but we repeat it here for the reader's sake.

\a {\bf Lemma 3.1.4.} {\it Let $k$ be an integer $\geq 6$ and $R$ be a 0-dimensional
scheme in  $\PPP (\Omega)$ such that $h^0({\cal E}_{k-5}\otimes {\cal
I}_{R})=0$. Let $A$ be the union in $\PPP (\Omega)$ of $R$ with $\pi ^{-1}(Y)$
where $Y$ denotes the union of $2k-4$ 2-fat points in $\PP 2$ supported at
general points; then, $h^0({\cal E}_{k+1}\otimes {\cal I}_{A})=0$.}

\a {\bf Proof.} Let $C$ be a smooth conic; we denote by $Z$ the scheme obtained
by $A$ specializing $k$ among the 2-fat points of $Y$ on $C$ and consider the
exact sequence:
$$0\rightarrow {\cal E}_{k-1}\otimes {\cal I}_{Res_{\pi ^{-1}C}Z}\rightarrow
{\cal E}_{k+1}\otimes {\cal I}_{Z}\rightarrow {\cal E}_{k+1}\otimes {\cal
I}_{Z\cap {\pi ^{-1}C},{\pi ^{-1}C}}\rightarrow 0;$$

\noindent since on $C$ there is a scheme of length $2k$, $H^0({\cal E}_{k+1}\otimes {\cal
I}_{Z\cap {\pi ^{-1}C},{\pi ^{-1}C}})=0$ so that $h^0({\cal E}_{k-1}\otimes
{\cal I}_{Res_{\pi ^{-1}C}Z})=h^0({\cal E}_{k+1}\otimes {\cal I}_{Z})\geq
h^0({\cal E}_{k+1}\otimes {\cal I}_{A})$.

\b The scheme $Res_{\pi ^{-1}C}Z$ is the union of $R$ and of the pull-back of
$k-4$ general 2-fat points and $k$ simple points on $C$. Now let $C'$ be
another smooth conic; we denote by $B$ the scheme obtained by $Res_{\pi
^{-1}C}Z$ specializing the $k-4$ 2-fat points on $C'$ and 4 among the simple
points on $C\cap  C'$; consider the exact sequence:
$$0\rightarrow {\cal E}_{k-3}\otimes {\cal I}_{Res_{\pi ^{-1}C'}B}\rightarrow
{\cal E}_{k-1}\otimes {\cal I}_{B}\rightarrow {\cal E}_{k-1}\otimes {\cal
I}_{B\cap {\pi ^{-1}C'},{\pi ^{-1}C'}}\rightarrow 0$$

\noindent and since on $C'$ there is a scheme of length $2k-4$, the third $H^0$ is 0 so
that $h^0({\cal E}_{k-3}\otimes {\cal I}_{Res_{\pi ^{-1}C'}B})=h^0({\cal
E}_{k-1}\otimes {\cal I}_{B})\geq h^0({\cal E}_{k-1}\otimes {\cal I}_{Res_{\pi
^{-1}C}Z})$.

\b The scheme $Res_{\pi ^{-1}C'}B$ is the union of $R$ and of the pull-back of
$k-4$ simple points on $C$ and $k-4$ simple points on $C'$. We denote by $D$
the scheme obtained by $Res_{\pi ^{-1}C'}B$ specializing the $k-4$ $C'$-points
on $C$ (for details, see \ref {I}); now consider the exact sequence:
$$0\rightarrow {\cal E}_{k-5}\otimes {\cal I}_{Res_{\pi ^{-1}C}D}\rightarrow
{\cal E}_{k-3}\otimes {\cal I}_{D}\rightarrow {\cal E}_{k-3}\otimes {\cal
I}_{D\cap {\pi ^{-1}C},{\pi ^{-1}C}}\rightarrow 0$$
and since on $C$ there is a scheme of length $2k-8$, the third $H^0$ is 0. Since
$Res_{\pi ^{-1}C}D$ is $R$, we finally get  $0= h^0({\cal E}_{k-5}\otimes {\cal
I}_{Res_{\pi ^{-1}C}D})=h^0({\cal E}_{k-3}\otimes {\cal I}_{D})\geq h^0({\cal
E}_{k-3}\otimes {\cal I}_{Res_{\pi ^{-1}C'}B})$, that is, $h^0({\cal
E}_{k+1}\otimes {\cal I}_{A})=0$. \qed

\a {\bf Lemma 3.1.5.} {\it If $d\leq {s\over 3}$ and $k\geq 10$, then $H(s,d,t,p,k)$
is true.}
\a {\bf Proof.}  Since $d\leq {s\over 3}$, by  Remark 2.2 the scheme $Z(s,d,t,p)$
specializes to $Z(s',0,t',p)$ where $s'=s-3d-6\big[ {s-3d\over 6} \big],\;
t'=t+d+\big[ {s-3d\over 6} \big]$.
\b Since $s'\leq 5$ and $p\leq \rho$, we have $10+\rho \geq 2s'+p \equiv \rho
(12)$ so that $2s'+p=\rho$, hence it is easy to see that the scheme
the union of $\tilde Y(s',0,0)$ and of $R_p$ specializes to
$R_{\rho}$ (see description of schemes $R_p$ in Notations 2.3). So finally the
scheme $Z(s',0,t',p)$ specializes to
$Z(0,0,u(k),r(k))$ and we conclude by semicontinuity that $H(s,d,t,p,k)$ holds.
\qed

\a {\bf Lemma 3.1.6.} {\it The statements $H(s,d,t,p,k)$ with $d> {s\over 3}$ are
true if both the statements $H(0,d,t,\rho(k),k)$ with any $d$ and the statements
$H(0,d,t,5,k)$ with $k\equiv 5$ {\rm (mod 6)} and with any $d$ are true.}

\a {\bf Proof.} Since $d>{s\over 3}$, Remark 2.2 allows us to say that the scheme
$Z(s,d,t,p)$ specializes to $Z(\sigma ,\delta ,\tau ,p)$ where $\sigma
=s-3\big[ {s\over 3} \big], \; \delta =d-\big[ {s\over 3} \big], \; \tau
=t+\big[ {s\over 3} \big]$; moreover, $0\leq \sigma \leq 2$ and $\delta \geq
1$.
Let $u=u(k)$, $\rho =\rho (k)$, and set $\delta =2e+j$, $0\leq j \leq 1$; we
have: $2\sigma +6\delta +12\tau +p=12(\tau +e)+6j+2\sigma +p=k(k+2)=12u+\rho$,
$\; 0\leq p\leq \rho$.
Since $10+\rho \geq 6j+2\sigma +p \equiv \rho$ (mod 12), we get
$6j+2\sigma+p=\rho$. Hence it is easy to see that if $(\sigma,j,p)\neq (0,1,5)$
the union of $\tilde Y(\sigma,j,0)$ and of $R_p$ specializes to
$R_{\rho}$, so that finally the scheme $Z(s,d,t,p)$ specializes to $Z(0,\delta
-j,\tau ,\rho)$, where $\delta -j \equiv 0$ (mod 2), and we conclude by semicontinuity. If $(\sigma,j,p)= (0,1,5)$
then $\rho=11$ so that $k\equiv 5$ (mod 6), and it is no longer true that
 the union of $\tilde Y(\sigma,j,0)$ and of $R_p$ specializes to
$R_{11}$, essentially because two double points in the plane do not specialize
to a triple point. \qed

\a {\bf Remark and notations 3.1.7.} Now our purpose is to prove the statements $H(0,d,t,\rho(k),k)$
with any $d$ and the statements $H(0,d,t,5,k)$ with $k\equiv 5$ (mod 6) and
any $d$, for $k\geq12$.
\b Notice that the number $d$ of double points in the statement $H(0,d,t,\rho(k),k)$ is necessarily even, since by assumption $6d+12t=k(k+2)-\rho (k)\equiv 0$ (mod 12) (see 3.1.2).
\b On the other hand, the number $d$ of double points in the statement $H(0,d,t,5,k)$ with $k\equiv 5$ (mod 6) is necessarily odd, since  (see 3.1.2 again) $k\equiv 5$ (mod 6) if and only if $\rho (k)=11$, and by assumption $6(d-1)+12t=k(k+2)-11\equiv 0$ (mod 12).

\par In the following we set $X(d,t,k):=Z(0,d,t,\rho) $ where
$6d+12t+\rho=k(k+2)$ and $k(k+2)=12u+\rho$, with $0\leq \rho \leq 11$; notice
that $t$ and $k$, as well as $d$ and $k$, determine $X(d,t,k)$, and $d$ is
always even. 
\par With $\bar R_{11}$ we denote in the following the generic union of the inverse image of a double
point with $R_5$. If $k\equiv 5$ (mod 6), we denote by $\bar X(d,t,k)$ the scheme
obtained by $X(d,t,k)$ substituting to $R_{11}$ the scheme $\bar R_{11}$.
\b So finally what we want to prove is that $X(d,t,k)$ is k-settled for any $k\geq12\;$
and that $\bar X(d,t,6l+5)$ is $(6l+5)$-settled for any $l\geq 2$.

\a { \maiu 3.2 Proof of the statement with no simple points.}

\a {\bf 3.2.1. Definition of standard step.} \par \noindent In the following $C$ denotes a smooth conic.
\par \noindent Let $R=\big ( x\pmatrix {1 \cr 2}+y\pmatrix {1 \cr 3}+ z\pmatrix {2
\cr 3}+ \pmatrix {e} \big )_{\, over\; C} +v\pmatrix {1 \cr 2\cr 3}+u\pmatrix {1
\cr 2}+R_\rho \;$ be a 0-dimensional scheme with length $R=h^0({\cal E}_{k+1})$
and length$(R\cap \pi ^{-1} C)\leq 2k$, and assume we want to prove that $R$ is
$k$-settled. We now define what a standard step is;
the idea is that we specialize (in the sense of a differential Horace
step) the maximum possible of triple points on
$C$, and if no more triple points are available, we specialize double points on
$C$, so that to get a scheme $\bar R$ with exactly $2k$ conditions on $C$; at
this point in order to prove that $R$ is $k$-settled it is enough to prove that
the residual scheme is $k-2$-settled:
\b {\bf standard step $\bf k \to k-2$}:
\b we ``add" over $C$ $\quad g\pmatrix {1 \cr 2\cr [3]}+ n\pmatrix {1 \cr [2]\cr
3}+ p\pmatrix {[1] \cr 2\cr
3}+r\pmatrix {1 \cr [2]}+s\pmatrix {[1] \cr 2}$ where:
\b $v\geq g+n+p$, $u\geq r+s$, $\;2x+3y+3z+e+3g+2n+p+2r+s=2k$,  $0\leq n+p\leq
1$, $0\leq s\leq 1$, and
finally $r=s=0$ if  it is possible to find $g,n,p$ such that
$2x+3y+3z+e+3g+2n+p=2k$.
\b The residue is:
$\big (\pmatrix {x+y+2z+r+2s}+g\pmatrix {1 \cr 2}+n\pmatrix {1 \cr 3}+ p\pmatrix
{2
\cr 3}\big)_{over \; C} +(v-g-n-p)\pmatrix {1 \cr 2\cr
3}+(u-r-s)\pmatrix {1 \cr 2}+R_\rho$.
\b Notice that this construction is possible if $\;2x+3y+3z+e \leq 2k$ and if
$v\geq g+n+p$, $u\geq r+s$.

\a {\bf Notation 3.2.2.} We recall here Definition  3.2 given in \ref {G-I}:

 Let $b,c,d,e,f,\rho ,k$ be integers $ \geq 0$; with
$Z(b,c,d,e,f,\rho ,k)$ we denote a $0$-dimensional subscheme of X, union of:
\b $b\pmatrix {1 \cr 2} + c(1) + d\pmatrix {1 \cr 3} + e{2 \choose 3} \quad {\rm
over}\quad C, \qquad{\rm and} \quad f \pmatrix {1 \cr 2\cr 3} + T_r \quad {\rm
general\quad over} \quad \PP 2, $

with the following assumptions:
\bigskip
\par \noindent
$\bf(0)_k$ \quad $2b+c+3d+3e \leq 2k$,  \quad$0 \leq d+e \leq 1$,
\par \noindent
$\bf(1)_k$ \quad $2(3b+c+4d+5e+6f)+r = k(k+2)$ \quad (i.e.,
$length(Z(b,c,d,e,f,r,k)) = h^0({\cal E}_{k+1}))$
\par \noindent
$\bf(2)_k$ \quad $\rho =0$  or $\rho =8$ if $k$ is even;  $\rho =3$  or $\rho
=11$ if $k$ is odd.

\a In \ref {G-I} it is proved that $Z(b,c,d,e,f,\rho ,k)$ is $k$-settled for $k\geq
12$: see Definition 3.2 and proof of Proposition
3.9 there.

\medskip \noindent If $k=6l+5$, i.e. if $\rho=11$, we denote here by  $\bar Z(b,c,d,e,f,11,k)$ the scheme
obtained by $Z(b,c,d,e,f,11,k)$ substituting to $R_{11}$ the scheme $\bar R_{11}$.

\a {\bf Lemma 3.2.3.} {\it Let $k=6l+5$. The schemes $\bar Z(b,c,d,e,f,11 ,k)$
 are k-settled for $k\geq 12$.}
\a {\bf Proof.} Use the proof of the fact that $Z(b,c,d,e,f,\rho ,k)$ is
$k$-settled given in \ref {G-I}, substituting to $R_{11}$ the scheme $\bar R_{11}$; it is hence enough to prove the initial
cases with this substitution. So it is enough to prove an analogous of Lemma
5.2 in \ref {G-I}, where 11 schemes of type  $Z(b,c,d,e,f,11,7)$  are proved to be
7-settled. We'll do the same here substituting to $R_{11}$ the scheme $\bar R_{11}$; so we want to prove that the schemes:
\b $\bar Z(0,9,0,1,2,11,7)$, $\bar Z(1,5,0,0,3,11,7)$, $\bar Z(1,6,0,1,2,11,7)$,
$\bar Z(1,7,1,0,2,11,7)$,
\b $\bar Z(2,2,0,0,3,11,7)$, $\bar Z(2,3,0,1,2,11,7)$, $\bar Z(2,4,1,0,2,11,7)$,
$\bar Z(3,5,0,0,2,11,7)$,
\b $\bar Z(4,2,0,0,2,11,7)$, $\bar Z(4,3,0,1,1,11,7)$, $\bar Z(1,1,1,0,3,11,7)$
\b are 7-settled. In all these cases we do 2 standard steps (see 3.2.1) and in all the 11 cases the
last residue is
$\big (\pmatrix {2}+\pmatrix {1 \cr 2}\big )_{\, over\; C} +\pmatrix {1/0 \cr 2}$ or
$\pmatrix
{5}_{\, over\; C} +\pmatrix {1/0 \cr 2}$ (where $\pmatrix {1/0 \cr 2}$ denotes the
scheme $R_5$). Now the second scheme specializes to the first one, since 5
points on a conic are general, so it is enough to prove that the first one is
1-settled.
\b We do a step $3 \to 1$: we ``add" over $C$ $\quad \pmatrix {1/0 \cr [2]}$ and
the
residue is $\pmatrix {1}+\pmatrix {1/0 }$ which is 1-settled (recall that the
scheme $\pmatrix {1/0 }_{\, over\; C}$ is just a point of $\PPP (\Omega)$). \qed

\a {\bf Lemma 3.2.4.} {\it If $d\leq k$, $X(d,t,k)$ is k-settled  for any $k\geq
12$  and $\bar X(d,t,k)$ is k-settled  for any $k\geq
12$  and $k\equiv 5$ (mod 6).}
\a {\bf Proof.} Let $u=u(k)$, $\rho =\rho (k)$. In \ref {G-I} it is proved that
$H^0({\cal E}_{k+1}\otimes
{\cal I}_{Z(b,0,0,0,f,\rho,k)})=0$ for $k\geq 12$, and Lemma 3.2.3 says that for
$k=6l+5$, the schemes $\bar Z(b,0,0,0,f,11,k)$ are k-settled for
$k\geq 12$. Now the scheme
$Z(b,0,0,0,f,\rho,k)$ is union in $\PPP (\Omega)$ of
$R_{\rho}$ with the inverse image of $f$ general triple points and of $b$ double
points whoose suports lie on a smooth conic, with $b\leq k$ (this is condition
$(0)_k$), so that we conclude by semicontinuity that $X(d,t,k)$ is k-settled
for any $k\geq 12$. Analougously for $\bar X(d,t,6l+5)$. \qed

\a {\bf Lemma 3.2.5.} {\it If $k\geq 16$ and $\; k< d\leq k+\big[{k-4\over 2}\big]$,
or if $k\geq 18$ and $\; k+\big[{k-4\over 2}\big]< d< 2k-4$, then $X(d,t,k)$ is
k-settled and if $k=6l+5$ $\bar X(d,t,k)$ is
$k$-settled.}

\a {\bf Proof.} We first prove the statement about $X=X(d,t,k)$. Let $C$ be a
smooth conic. We do an Horace step. We specialize
on $C$ $k$ double points; now on
$C$ there is a scheme of length $2k$, so that $H^0({\cal E}_{k+1}\otimes {\cal
I}_{X\cap {\pi ^{-1}C},{\pi
^{-1}C}})=0$, hence $h^0({\cal E}_{k-1}\otimes {\cal I}_{Res_{\pi
^{-1}C}X})=h^0({\cal E}_{k+1}\otimes {\cal I}_{X})$. The residual scheme $Res_{\pi ^{-1}C}X$ is the generic union of $\tilde Y(0,d-k,t)$ with the pull back of $k$
points on $C$ and
with $R_{\rho}$.
\b Now we do an HD step; this time we need to have $2k-4$ points of $\PP 2$ in
total on $C$, so we still need $k-4$.

\smallskip \noindent {\it case 1: $d - k\leq \big[{k-4\over 2}\big]$}:
\b we add on $C$ $(d -k) \pmatrix {1 \cr [2]}+
  g\pmatrix {1 \cr 2\cr [3]}+ h\pmatrix {1 \cr [2]\cr 3}+ i\pmatrix {[1] \cr
2\cr 3}$ where $k+2(d-k)+3g+2h+i=2k-4$ and $0\leq h+i\leq 1$. The HD
residue is the generic union of
$\tilde Y(0,0,t-g-h-i)$ of $R_{\rho}$ and of $(d -k) (1)+g\pmatrix
{1 \cr 2}+ h\pmatrix {1 \cr 3}+ i\pmatrix { 2\cr 3}$ on $C$, which in the
notation of \ref {G-I} is  a generalization of $Z(g,d-k,h,i,t-g-h-i,\rho,k-4)$; in
fact, it is easy to check that
conditions $(1)_{k-4}$ and $(2)_{k-4}$ are automatically verified, while the condition 
$(0)_{k-4}$ is true for
$k\geq 4$. Moreover, $t-g-h-i \geq 0$ if $k\geq 16$.

\smallskip \noindent {\it case 2: $d - k> \big[{k-4\over 2}\big]$}:
\b we set $k-4=2m+l$, $l=0,1$, and we add on $C$ $m \pmatrix {1 \cr [2]}+
 l\pmatrix {[1] \cr 2}$ where $k+2m+l=2k-4$. The HD residue the generic union of
$\tilde Y(0,d-k-m-l,t)$, of 
$R_{\rho}$ and of $m(1)+l(2)$ on $C$.
\b Now we do another HD step: we need $2k-8$ points of $\PP 2$ in
total on $C$, so we add on $C$ $(d -k-m-l) \pmatrix {1 \cr [2]}+
  g\pmatrix {1 \cr 2\cr [3]}+ h\pmatrix {1 \cr [2]\cr 3}+ i\pmatrix {[1] \cr
2\cr 3}$ where $m+2l+2(d-k-m-l)+3g+2h+i=2k-8$ and $0\leq h+i\leq 1$. The HD
residue is the generic union of
$\tilde Y(0,0,t-g-h-i)$, of 
$R_{\rho}$ and of $(d -k-m-l)
(1)+g\pmatrix
{1 \cr 2}+ h\pmatrix {1 \cr 3}+ i\pmatrix { 2\cr 3}$ on $C$, which in the
notation of \ref {G-I} is  a generalization of
$Z(g,d-k-m-l,h,i,t-g-h-i,\rho,k-6)$; in fact,
conditions $(1)_{k-6}$ and $(2)_{k-6}$ are automatically verified, while
$(0)_{k-6}$ is true for
$k\geq 16$. Moreover, it is easy to check that $t-g-h-i\geq 0$ if $k\geq 18$.

\smallskip The conclusion follows by \ref {G-I} , where, as said above, it is proved
that these schemes $Z(b,c,d,e,f,r,k)$ are $k$-settled for $k\geq 12$.
\b The statement about $\bar X(d,t,6l+5)$  is proved exactly in the same way
using Lemma 3.2.3 instead of \ref {G-I}. \qed

\a {\bf 3.2.6.} In the following Lemma 3.2.7 we treat the case of double and
triple points with a lot of double points. Hence it is convenient to give some
definitions.

\a Let $k,\, n$ be integers with  $1\leq n\leq {k\over 6}$ and let $R$ be a
0-dimensional scheme in  $\PPP (\Omega)$ such that $h^0({\cal
E}_{k+1-6n}\otimes {\cal I}_{R})=0$. Let $A$ be the union in $\PPP (\Omega)$ of
$R$ with the inverse image of the union of $(2k-4)+(2(k-6)-4)\dots
+(2(k-6(n-1))-4)=n(2k-6n+2)$ double points in $\PP 2$ supported at general
points; then, Lemma 3.1.4 applied $n$ times gives $h^0({\cal E}_{k+1}\otimes {\cal
I}_{A})=0$ (the condition $n\leq {k\over 6}$ assures that $k+1-6n \geq 1$, and
also that each addend in the sum is positive, since $n\leq {k+4\over 6}$).
\b For all $k \geq 6$ and $1\leq n\leq {k+4\over 6}$ we set $$\; \alpha
(n,k):=\sum _{ i=0}^{n-1}\big( 2(k-6i)-4\big) =n(2k-6n+2).$$
\b If we fix a $k \geq 6$ the function $\; \alpha (n,k)$ is hence increasing as
long as it is defined, and strictly increasing if $n< {k+4\over 6}$.
\b Now consider a scheme $X(d,t,k)$ with $k \geq 6$; we can set

$$ \bar n= n(d,k):=max \{n, \; 1\leq n\leq {k\over 6},\;  d\geq \alpha (n,k)\}\;
{\rm if} \; d\geq 2k-4,$$
$$ \bar n= n(d,k):=0 \; {\rm if} \; d<2k-4. $$

\b Let $d\geq 2k-4$ and let $m$ be an integer, $1\leq m\leq \bar n$; we have
seen above that
$H(0,d,t,\rho,k)$ is true if $H(0,d-\alpha(m,k),t,\rho,k-6m)$ is true. Moreover,
the scheme $X(d-\alpha(m,k),t,k-6m)$ verifies

$$\alpha(\bar n +1 -m, k-6m)>d-\alpha(m,k)\geq \alpha(\bar n -m, k-6m)\quad
\quad (*)$$

 so that $$ n(d-\alpha(m,k),k-6m)= n(d,k) -m\quad \quad (**).$$
\b In fact, one has $\alpha (\bar n,k)-\alpha (m,k)=\alpha(\bar n -m, k-6m)$, so
that the second inequality is clear. For the first one, there are two
possibilities:
\b i) $\alpha(\bar n +1, k)$ is defined and $>\alpha(\bar n, k)$, i.e. $\bar
n+1<
{k+4\over 6}$; then  by definition $ \alpha (\bar n+1,k)>d$ so that $\alpha(\bar
n +1-m, k-6m)>d-\alpha(m,k)$.
\b ii) $\alpha(\bar n +1, k)$ is not defined, i.e. $\bar n+1> {k+4\over 6}$, or
$\alpha(\bar n +1, k)=\alpha(\bar n , k)$, i.e. $\bar n+1= {k+4\over 6}$. 
\b Since $\bar n\leq {k\over 6}$ we have: $\bar n= {k-1\over 6}$ or $\bar n= {k\over 6}$
in the first case, $\bar n= {k-2\over 6}$ in the second. We recall that
$6d+12t+\rho (k)=k(k+2)$, hence $d\leq {k(k+2)-\rho \over 6}=2u(k)$.
\b If $\bar n= {k-2\over 6}$ then $k\equiv 2$ (mod 6) so that $\rho=8$; hence,
$\alpha (\bar n,k)={(k-2)(k+4)\over 6}=2u(k)$.
\b If $\bar n= {k-1\over 6}$ then $k\equiv 1$ (mod 6) so that $\rho=3$; hence,
$\alpha (\bar n,k)={(k-1)(k+3)\over 6}=2u(k)$.
\b If $\bar n= {k\over 6}$ then $k\equiv 0$ (mod 6) so that $\rho=0$; hence,
$\alpha (\bar n,k)={k(k+2)\over 6}=2u(k)$.
Since by definition $d\geq \alpha (\bar n,k)$, in each of the three cases we get
$d=\alpha (\bar n,k)$; hence the first inequality in $(*)$ becames  $\alpha(\bar
n +1 -m, k-6m)>\alpha(\bar n-m, k-6m)$ which is true (notice that $\bar n +1
-m\leq \bar n$ so that
$\alpha(\bar n +1 -m, k-6m)$ is defined).

\a {\bf Lemma 3.2.7.} {\it If $k\geq 12$ then
$X(d,t,k)$ is k-settled and when $k=6l+5$ also  $\bar X(d,t,k)$ is k-settled.}

\a {\bf Proof.}  We prove the statement about $X(d,t,k)$, since the statement
about $\bar X(d,t,k)$ can be proved exactly in the same way.
\b If $\; d< 2k-4$, then $X(d,t,k)$ is $k$-settled for $k\geq 18$ by Lemma 3.2.4
and Lemma 3.2.5, and for $12\leq k \leq 17$ by Lemma 3.3.1.
\b If $\; d\geq 2k-4$, write $k=6h+j$,
$0\leq j \leq 5$, and let $\bar n = n(d,k)$. If $\bar n\leq h-3$, then
$k-6\bar n\geq 18+j$, so that $H(0,d-\alpha(\bar n,k),t,\rho,k-6\bar n)$ is true
by Lemma 3.2.4
and Lemma 3.2.5, since by $(*)$  $d-\alpha(\bar n,k)<2(k-6\bar n)-4$. If $\bar n\geq h-2$,
then we
apply  $h-2$
times Lemma 3.1.4, and in this way we see that $H(0,d,t,\rho,k)$ is true if
$H(0,d-\alpha(h-2,k),t,\rho,12+j)$ is true. So we conclude by Lemma 3.3.1.

\a { \maiu 3.3 Initial cases and proof of of Theorem 1.1.}

\a {\bf Lemma 3.3.1.} {\it If $k=12+j$ with $0\leq j \leq 5$ then $X(d,t,k)$ is
k-settled and when $k=17$ also  $\bar X(d,t,k)$ is k-settled.}

\a {\bf Proof.} If $d\leq k$ the statement is true by Lemma 3.2.4. If $k=16,17$ and
$k<d\leq k+\big[{k-4\over 2}\big]$, the statement is true by Lemma 3.2.5. We are
hence
going to prove that the union $X(d,t,k)$ of $d$ double and $t$
triple points and of $R_\rho$ verifies $H^0({\cal E}_{k+1}\otimes {\cal
I}_{X(d,t,k)})=0$, where $d+2t={k(k+2)- \rho \over 6}=2u(k)$ so
that $d$ is even, and $d> k$ if $12\leq k \leq 15$, $d>k+\big[{k-4\over 2}\big]$
if $k=16,17$. For $k=17$ we are also going to prove that the same holds
substituting to $R_{11}$ the scheme $\bar R_{11}$.
When $j=0,1,3,4$ the case $t=0$ is
proved in \ref {I}, section 2 so that we'll assume $t>0$ for $k=12,13,15,16$.
Recall that at each step $l\to l-2$, which can be an Horace or a differential Horace 
step, the divisor is the pull back of a smoth conic $C$ and the
points of $\PP 2$ needed on $C$ are $2l$. All the assertions about unions of double points plus a scheme $R_{\rho}$ are proved in \ref {I} section 2.

\a Assume $0\leq t\leq \big [{2k \over 3
}\big ]$. We define the following algorithm {\bf (A)(t,k)} applying one standard
step (see 3.2.1) to $X(d,t,k)$, and then another standard
step to the residue (in these assumptions the standard steps are particularly simple):
\b  $step \; \;k \to k-2$: We ``add" over $C$ $\quad t\pmatrix {1 \cr 2\cr
[3]}+b\pmatrix {1 \cr [2]}+ a\pmatrix {[1] \cr 2}$ where
$3t+2b+a=2k$,  $a=0,1$; since the condition $3t\leq 2k$ is verified by
assumption, this is possible if $d-a-b\geq 0$ is true.
\b The residue is $\big (t\pmatrix {1 \cr 2}+\pmatrix {b+2a}\big
)_{\, over\; C} +(d-a-b)\pmatrix {1 \cr 2}+R_\rho$.

\a $step \; \;k-2 \to k-4$: We ``add" over $C$ $\quad h\pmatrix {1 \cr [2]}+
i\pmatrix {[1] \cr 2}$ where
$2t+b+2a+2h+i=2(k-2)$,  $i=0,1$; this is possible if the conditions $2t+b+2a\leq
2(k-2)$ and $d-a-b-h-i\geq 0$ are verified.
\b The residue is $\pmatrix {t+h+2i}_{\, over\; C}
+(d-a-b-h-i)\pmatrix {1 \cr 2}+R_\rho$.
We set $w=w(t,k):=t+h+2i$, and $q=q(t,k):=d-a-b-h-i$; we are reduced to prove
that the scheme $\pmatrix {w}_{\, over\; C}+q\pmatrix {1 \cr 2}+R_\rho$ is
$(k-4)-settled$.

\a  Assume $ t> \big [{2k \over 3 }\big ]$. We
define the following algorithm {\bf (C)(t,k)} applying one standard
step to $X(d,t,k)$ (again in these assumptions the standard step is very simple):

\b step $k \to k-2$:
\b we ``add" over $C$ $\quad g\pmatrix {1 \cr 2\cr [3]}+ n\pmatrix {1 \cr [2]\cr
3}+ p\pmatrix {[1] \cr 2\cr
3}$ where $3g+2n+p=2k$,  $n+p=0,1$ (this is always possible).
\b The residue is $\big (g\pmatrix {1 \cr 2}+n\pmatrix {1 \cr 3}+ p\pmatrix {2
\cr 3}
\big )_{\, over\; C} +(t-g-n-p)\pmatrix {1 \cr 2\cr 3}+d\pmatrix {1 \cr 2}+R_\rho$.

\a In the following we apply the standard step 3.2.1, (A)(t,k) and (C)(t,k) a number
of times; easy calculations
assure that the conditions respectively $\;2x+3y+3z+e \leq 2k$, $v\geq g+n+p$,
$u\geq r+s$ for the standard step, and $2t+b+2a\leq 2(k-2)$, $d-a-b-h-i\geq 0$
for (A)(t,k) are verified (while (C)(t,k) is always possible).

\a ${\bf case\; k=12}$: here $\rho=0$; we have to treat the cases $1\leq t \leq
7$.
\b For $1\leq t \leq 7$ (A)(t,12) gives
$(w,q)=(10,10)$ for $t=6,7$, $(w,q)=(7,11)$
for $2\leq t\leq 5$, $(w,q)=(4,12)$ for $t=1$.

\b In all three cases we do a standard step $8 \to 6$, and the residue is
$\pmatrix {3}_{\, over\; C}+7\pmatrix {1 \cr 2}$, which specializes
to $8$ double points and it is 6-settled, or $\pmatrix
{6} _{\, over \; C}+6\pmatrix {1\cr 2}$.
In the last case, we do another standard step $6\to 4$, and the residue is
$\pmatrix {3} _{\, over\;
C}+3\pmatrix {1\cr 2}$, which specializes to $4$ double points and it is
4-settled.

\a ${\bf case\; k=14}$: here $\rho=8$; we have to treat the cases $0\leq t \leq
10$.

\b If $t\leq 9$ we use the first step of (A)(t,14), so that we have to prove
that
the residue $\big (t\pmatrix {1 \cr 2}+\pmatrix {b+2a}\big
)_{\, over\; C} +(d-a-b)\pmatrix {1 \cr 2}+R_8$ is 12-settled; recall that $R_8$ is
a scheme $\pmatrix {1 \cr 3}$.
\b step $12 \to 10$: We ``add" over $C$ $\quad h'\pmatrix {1 \cr [2]}+
i'\pmatrix {[1] \cr 2}+\pmatrix {1 \cr [3]}$ where
$2t+b+2a+2h'+i'+3=2(k-2)=24$,  $i'=0,1$; this is possible since the conditions
$2t+b+2a+3\leq
24$ and $d-a-b-h'-i'\geq 0$ are verified.
\b The residue is $\pmatrix {t+h'+2i'+1}_{\, over\; C}
+(d-a-b-h'-i')\pmatrix {1 \cr 2}$.
We set $w'=w'(t):=t+h'+2i'+1$, and $q'=q'(t):=d-a-b-h'-i'$; we are reduced to
prove
that the scheme $\pmatrix {w'}_{\, over\; C}+q'\pmatrix {1 \cr 2}$ is
$10-settled$; an easy calculation shows that $\quad \quad (w',q')=(12,16)\;$ for
$\;8\leq t \leq 9$,
$\;(w',q')=(9,17)$ for $\;4\leq t\leq 7$, $\; (w',q')=(6,18)$ for $\;0\leq
t\leq 3$.
These three configurations specialize to the first residual scheme of (A)(t,12)
(i.e.  $\big (t\pmatrix {1 \cr 2}+\pmatrix {b+2a}\big )_{\, over\; C}
+(d-a-b)\pmatrix {1 \cr 2}$) obtained respectively in cases $t=1,2,4$, and
we have proved that they are 10-settled.

\b If $t=10$ apply (C)(10,14); it is now enough to prove that $\big (9\pmatrix
{1 \cr 2}+\pmatrix {2 \cr 3}\big )_{\, over\; C} +16\pmatrix {1 \cr 2}+R_8$ is
12-settled.
\b step $12 \to 10$: we ``add" over $C$ $\quad \pmatrix {1 \cr [3]}$ (which is
the scheme $R_8$) and the residue is
 $\pmatrix {12}_{\, over\; C}+16\pmatrix {1 \cr 2}$, which is the case
$(w',q')=(12,16)\;$ previously treated.

\a ${\bf case\; k=16}$: here $\rho=0$; we have to treat the cases $1\leq t \leq
12$.
\b For $1\leq t \leq 10\;$ (A)(t,16) gives
$\quad (w,q)=(12,24)\;$ for $\;6\leq t \leq 10$,
$\;(w,q)=(9,25)$ for $\;2\leq t\leq 5$, $\;
(w,q)=(6,26)$ for $t=1$. Now we do a standard step $12 \to 10$, and if $w=12$,
the residue is
 $\pmatrix {6}_{\, over\; C}+18\pmatrix {1 \cr 2}$, while if $w=9$ or $w=6$, the
residue is $\pmatrix {9}_{\, over\; C}+17\pmatrix {1 \cr 2}$, and these are the
cases $(w',q')=(9,17)$ or
$(6,18)$ previously treated in $k=14$.

\b For $11\leq t \leq 12\;$ apply (C)(t,16); it is now enough to prove that
$\big (10\pmatrix {1 \cr 2}+\pmatrix {1 \cr 3}\big )_{\, over\; C}+(t-11)\pmatrix {1 \cr
2 \cr 3} +d\pmatrix {1 \cr 2}$ is 14-settled. We now do two standard steps more,
$14\to 12$ and $12\to 10$, and the last residue is in both cases $\pmatrix {6}
_{\, over \; C}+18\pmatrix {1 \cr 2}$, which we have just recalled is 10-settled.

\a ${\bf case\; k=13}$: here $\rho=3$; we have to treat the cases $1\leq t \leq
9$.
\b For $1\leq t \leq 8$ (A)(t,13) gives
$(w,q)=(12,12)$ for $t=8$, $(w,q)=(9,13)$
for $4\leq t\leq 7$, $(w,q)=(6,14)$ for $1\leq t\leq 3$.
\b In all three cases we do a standard step $9 \to 7$, and if $w=12$, the
residue is
 $\pmatrix {3}_{\, over\; C}+9\pmatrix {1 \cr 2}+R_3$, which specializes
to $10$ double points +$R_3$ and it is 7-settled.
\b If $w=9$ or $w=6$, the residue is $\pmatrix {6}_{\, over\; C}+8\pmatrix {1 \cr
2}+R_3$; we do another standard step $7\to 5$, and the residue is $\pmatrix {4}
_{\, over\;
C}+4\pmatrix {1\cr 2}+R_3$, which specializes to $5$ double points +$R_5$ and it
is 5-settled.
\a For $t=9$ apply (C)(9,13); it is now enough to prove that $\big (8\pmatrix {1
\cr 2}+\pmatrix {1 \cr 3}\big )_{\, over\; C} +14\pmatrix {1 \cr 2}+R_3$ is 11-settled.
We
do a standard step $11 \to 9$, and the residue is
 $\pmatrix {12}_{\, over\; C}+12\pmatrix {1 \cr 2}+R_3$, which is one of the
previous cases.

\a ${\bf case\; k=15}$: here $\rho=3$; we have to treat the cases $1\leq t \leq
13$.
\b For $1\leq t \leq 10\;$ (A)(t,15) gives
$\; (w,q)=(13,19)\;$ for $\;8\leq t \leq 10$,
$\;(w,q)=(10,20)$
for $\;4\leq t\leq 7$, $\; (w,q)=(7,21)\;$ for $1\leq t\leq 3$.
These three configurations specialize to the first residual scheme of (A)(t,13)
(i.e.  $\big (t\pmatrix {1 \cr 2}+\pmatrix {b+2a}\big )_{\, over\; C}
+(d-a-b)\pmatrix {1 \cr 2}+R_3$) obtained respectively in cases $t=1,2,4$, and
we have proved that they are 11-settled.
\a For $11 \leq t \leq 13\;$ apply (C)(t,15); it is now enough to prove that
$\big (10\pmatrix {1 \cr 2}\big )_{\, over\; C}+(t-10)\pmatrix {1 \cr 2 \cr 3}
+d\pmatrix {1 \cr 2}+R_3$ is 13-settled. If $11 \leq t \leq 12\;$, this scheme
specializes to the  scheme used in the first step of (A)(t,13) (i.e.  $\big
(t\pmatrix {1 \cr 2 \cr 3}+b\pmatrix {1 \cr 2}+ a\pmatrix {1 \cr
2}\big )_{\, over\; C}
+(d-a-b)\pmatrix {1 \cr 2}+R_3$) respectively in
cases $t=1$ (where $b=11$, $a=1$) and $t=2$ (where $b=10$, $a=0$), and we have
proved that they are 13-settled.
If $t=13$ we do three standard steps from 13 to 7 and the last
residue is $\pmatrix {6}_{\, over\; C}+8\pmatrix {1 \cr 2}+R_3$; in $k=13$ we have
proved that it is 7-settled.

\a ${\bf case\; k=17}$: here $\rho=11$; we have to treat the cases $0\leq t \leq
14$.
\b We first prove that $X(d,t,17)$ is 17-settled.
\b For $0\leq t \leq 11\;$ we have $d\geq 30$, and by Lemma 3.1.4 $X(d,t,17)$ is
17-settled if $X(d-30,t,11)$ is 11-settled. If $t=11$, $d=30$, so that $X(d-30,t,11)$ is
11-settled by \ref {G-I}.
\b If $0\leq t \leq 7\;$ (A)(t,11) applied to the scheme $X(d-30,t,11)$ gives
$\quad (w,q)=(8,6)\;$ for $\;4\leq t \leq 7$, $\quad (w,q)=(5,7)\;$ for $\;0\leq
t \leq 3$.
Now we do another standard step, and the residue is $\pmatrix {3}_{\, over\;
C}+3\pmatrix {1 \cr 2} +R_{11}$, resp. $\pmatrix {6}_{\, over\; C}+2\pmatrix {1 \cr
2} +R_{11}$.
\b step $5 \to 3$: we ``add" over $C$ $\quad r\pmatrix {1 \cr [2]}+s\pmatrix
{[1] \cr 2}+\pmatrix {1/0 \cr 2 \cr [3]}$ with $0\leq s\leq 1$ so to have 10
conditions on $C$, and
in both cases the residue is $\big (\pmatrix {2}+\pmatrix {1/0 \cr 2}\big )_{\, over
\; C}+\pmatrix {1 \cr 2}$.
\b step $3 \to 1$: we ``add" over $C$ $\quad \pmatrix {1 \cr [2]}$ and the
residue is $\pmatrix {1}+\pmatrix {1/0 }$ which is 1-settled (recall that the
scheme $\pmatrix {1/0 }_{\, over\; C}$ is just a point of $\PPP (\Omega)$).

\medskip \b For $8\leq t \leq 10$, we do 3 standard step from 11 to 5, and the
residue is $\pmatrix {3}_{\, over\; C}+3\pmatrix {1 \cr 2} +R_{11}$ if $t=8,9$,
which has been treated above, or $\big (\pmatrix {2}+\pmatrix {1 \cr 3} \big
)_{\, over\; C}+2\pmatrix {1 \cr 2}+R_{11}$ if $t=10$, in which case we proceed as
follows:
\b step $5 \to 3$: we ``add" over $C$ $\quad \pmatrix {1/0 \cr 2 \cr
[3]}+\pmatrix {1 \cr [2]}$ and the
residue is $\big (\pmatrix {2}+\pmatrix {1/0 \cr 2} \big )_{\, over\; C}+\pmatrix {1
\cr
2}$ treated above.

\b If $12\leq t \leq 14\;$, we do 3 standard steps $17\to 15$, $15\to 13$ and
$13\to 11$, and the
residue is $\pmatrix {6}_{\, over\; C}+20\pmatrix {1 \cr 2} +R_{11}$ in all the
three
cases; we go on with 3 other standard steps $11\to 9$, $9\to 7$ and $7\to 5$,
and the residue is
$\pmatrix {6}_{\, over\; C}+2\pmatrix {1 \cr 2} +R_{11}$, already treated above.
\b  Now we want to prove that $\bar X(d,t,17)$ is 17-settled, so we substitute
to $R_{11}$ the scheme $\bar R_{11}$, and, if $t\neq 11$, we do the same steps from 17 to 5
as above; now we have to prove that the following schemes are 5-settled:
\b $\pmatrix {6}_{\, over\; C}+2\pmatrix {1 \cr 2} +R_5 +\pmatrix {1 \cr 2}$;
\b $\pmatrix {3}_{\, over\; C}+3\pmatrix {1 \cr 2} +R_5 +\pmatrix {1 \cr 2}$;
\b $\big (\pmatrix {2}+\pmatrix {1 \cr 3} \big )_{\, over\; C}+2\pmatrix {1 \cr
2}+R_5 +\pmatrix {1 \cr 2}$.
We do a standard step $5\to 3$ and the
residue is $\pmatrix {2}_{\, over\; C}+\pmatrix {1/0 \cr 2}+\pmatrix {1 \cr 2}$ in
the first case, $\pmatrix
{5}_{\, over\; C} +\pmatrix {1/0 \cr 2}$ in the second and third case (here
$\pmatrix {1/0 \cr 2}$ denotes the scheme $R_5$); in Lemma 3.2.3 we proved that
both are 3-settled. 
\b Now let $t=11$. By Lemma 3.1.4 it is enough to prove that $\bar X(0,11,11)$ is 11-settled. In order to prove this, we do 4 standard steps from 11 to 3 and the residue is $\big (\pmatrix {2}+\pmatrix {1 \cr 2} \big )_{\, over\; C} +\pmatrix {1/0 \cr 2}$ which is proved to be 3-settled in Lemma 3.2.3.
\qed

\a {\bf Corollary 3.3.2.} {\it $H(s,d,t,p,k)$ is true for $k\geq 12$}.

\a {\bf proof.} It follows by Lemma 3.1.5, Lemma 3.1.6 and Lemma 3.2.7. \qed

\a {\bf Proof of Theorem 1.1.} First let us check that in our assumptions $Y(a,b,c)$ has maximal Hilbert function (mHf for short). Notice that $Y(a,b,c)$ has mHf if and only if $h^0({\cal I}_{Y(a,b,c)}(v-1))=0$ and $h^0({\cal I}_{Y(a,b)}(v))={v+2\choose 2}-l$ where $l=l(Y(a,b,c))=a+3b+6c$ and $v=v(Y(a,b,c))$ (see 2.1).
\par By \ref{M}, a general fat point scheme $Z=m_1P_1+\dots+m_rP_r$, $4\geq m_1\geq\dots\geq m_r\geq0$, has mHf in any degree $k$ such that $k\geq m_1+m_2+m_3$. Since for our schemes $Y(a,b,c)$ one always has $m_1+m_2+m_3\leq 9$, we get that $Y(a,b,c)$ has mHf for $v\geq 10$, i.e. for $l>55$. 

\medskip Let $k$ be the integer such that $(k-1)(k+1) <2a+6b+12c\leq k(k+2)$, so that $w(Y(a,b,c))=k$. Rephrasing what is done in the proof of Lemma 2.6 for the analogous statements with $c=0$, it is easy to show that 
\b i) there exists $(s,d,t,p)\in \Lambda _{k-1}$ such that $Z(s,d,t,p)\subseteq \tilde Y(a,b,c)$;
\b ii) there exists $(s',d',t',p')\in \Lambda _{k}$ such that $Z(s',d',t',p')\supseteq \tilde Y(a,b,c)$.

\medskip  If $2l> 12\cdot 14$, i.e. $l>84$, then $k\geq 13$, hence $(s,d,t,p)\in \Lambda _{k-1}$, respectively $(s',d',t',p')\in \Lambda _{k}$, implies that $H(s,d,t,p,k-1)$, respectively that $H(s',d',t',p',k)$ is true, i.e. $H^0({\cal E}_{k}\otimes {\cal I}_{Z(s,d,t,p)})=0$, respectively $H^0({\cal E}_{k+1}\otimes {\cal I}_{Z(s',d',t',p')})= 0$ (see Corollary 3.3.2 and 3.1.2). 
\par So we see, exactly as in the proof of Lemma 2.6, that $\mu_{k-1}(Y(a,b,c))$ is injective and $\mu_{k}(Y(a,b,c))$ is surjective and we conclude that $Y(a,b,c)$ is minimally generated (see 2.1). 
If $k=12$ (i.e. $143<2l\leq 168)$ and $v=12$ (i.e. $78<l\leq 91)$, then in order to prove that $Y(a,b,c)$ is minimally generated it is enough to prove that $\mu_{12}$ is surjective (see 2.1), and this is true again by ii) and Corollary 3.3.2; hence $Y(a,b,c)$ has the expected resolution also for $79\leq l\leq 84$. \qed

\bigskip 
\bigskip 
\noindent {\bf References}

\bigskip
\ref {A-H} J. Alexander and A. Hirschowitz, {\it An asymptotic vanishing
theorem for generic unions of multiple points}, Invent.
Math. 140 (2000), no. 2, 303--325.

\ref {B} E. Ballico, {\it Generators for the homogeneous ideal of $s$
general points of $\Bbb {P}^3$}, J. Algebra 106 (1987),
no. 1, 46--52.

\ref {B-G} E. Ballico and A. V. Geramita, {\it The minimal free
resolution of the ideal of $s$ general points in $\Bbb {P}^3$},
in: Proceedings of the 1984 Vancouver Conference in Algebraic
Geometry, CMS Conference Proceedings Vol. 6, 1--11,
Canadian Math. Soc. and Amer. Math. Soc, Providence, RI, 1986.

\ref {Ca} M.V. Catalisano, {\it "Fat" points on a conic},
Comm. Algebra 19 (1991), 2153--2168.

\ref {C} K. A. Chandler, {\it A brief proof of a maximal rank theorem
for generic double points in projective space}, Trans.
Amer. Math. Soc. 353 (2000), no. 5, 1907--1920.

\ref {C-C-M-O} C. Ciliberto, F. Cioffi, R. Miranda and F. Orecchia.
{\it Bivariate Hermite interpolation and linear systems of
plane curves with  base fat points}, in: Computer mathematics,
87--102, Lecture Notes Ser. Comput., 10, World Sci. Publishing,
River Edge, NJ,  2003.

\ref {E} L. Evain, {\it Calculs de dimensions de syst\`emes lin\'eaires de courbes planes par collisions de gros points}, C.R.Acad.Sci.Paris, t.325, S\'erie I (1997), 1305--1308.

\ref {E2} L. Evain, {\it Computing limit linear series with infinitesimal
methods}, preprint 2004 (arXiv:math. AG/0407143).

\ref {F} S. Fitchett, {\it  Maps of linear systems on blow ups
of the Projective Plane}, J. Pure Appl. Algebra 156 (2001), 1--14.

\ref {F-H-H} S. Fitchett, B. Harbourne and S. Holay, {\it Resolutions of Fat Point Ideals Involving Eight
General Points of $\PP2$}, J. Algebra 244 (2001), 684--705.

\ref {G-I} A. Gimigliano and M. Id\`{a}, {\it The ideal resolution for
generic $3$-fat points in $\Bbb {P}^2$}, J. Pure Appl.
Algebra 187 (2004), no. 1-3, 99--128.

\ref {G-M} A. V. Geramita and P. Maroscia, {\it The ideal of forms
vanishing at a a finite set of points in $\Bbb {P}^n$},
J. Algebra 90 (1984), no. 2, 528--555.

\ref {Ha1} B.\ Harbourne, {\it The Ideal Generation
Problem for Fat Points}, J.\ Pure Appl.\ Alg.\ 145(2), 165--182 (2000).

\ref {Ha2} B. Harbourne, {\it Anticanonical rational surfaces},
Trans. Amer. Math. Soc. 349, 1191--1208 (1997).

\ref {Ha3} B. Harbourne, {\it An Algorithm for Fat Points
on $\PP2$}, Can. J. Math. 52 (2000), 123--140.

\ref {H-H-F} B.\ Harbourne, S.\ Holay and S.\ Fitchett,
{\it Resolutions of ideals of
quasiuniform fat point subschemes of $\Bbb {P}^2$},
Trans.\ Amer.\ Math.\ Soc.\ 355 (2003), no.\ 2, 593--608.

\ref {H-R} B. Harbourne and J. Ro\'e. {\it Linear systems with multiple base points in $\PP 2$},
Adv.\ Geom. 4 (2004), 41--59.

\ref {Hi} A. Hirschowitz, {\it La m\'{e}thod d'Horace pour
l'interpolation \`{a} plusieurs variables}, Manuscripta Math. 50
(1985),
337--388.

\ref {H-S}  A.Hirschowitz - C.Simpson:  {\it La r\'esolution minimale 
de l'id\'eal d'un arrangement g\'en\'eral d'un grand nombre 
de points dans ${\Bbb P}^n$}, Invent. Math. 126 (1996), 467--503. 

\ref {I} M. Id\`{a}, {\it The minimal free resolution for the first
infinitesimal neighborhood of $n$ general points in the plane}, J. Algebra 216 (1999), no. 2, 741--753.

\ref {I2} M.Id\`a: {\it Generators for the generic rational space 
curve: low degree cases}. "Commutative Algebra and Algebraic 
Geometry, Proc. of the 
Ferrara Conference in honour of M.Fiorentini, 1997", Lecture Notes in 
Pure and Applied Math. Dekker  {\bf 206} (1999), 169--210.

\ref {M} T. Mignon, {\it Syst\'{e}mes de courbes plane \`{a}
singularit\'{e}s impos\'{e}es: le cas des multiplicit\'{e}s
inf\'{e}rieures
ou \'{e}gales \`{a} quatre}, J. Pure Appl. Algebra 151 (2000), no. 2, 173--195.

\ref {N} M. Nagata, {\it On rational surfaces, II},
Mem.\ Coll.\ Sci.\
Univ.\ Kyoto, Ser.\ A Math.\ 33 (1960), 271--293.

\ref {R} J. Ro\'{e}, {\it Limit linear systems and applications},
e-print arXiv:math.AG/0602213.

\ref {Y} S. Yang, {\it Linear series in $\PP2$ with base points of
bounded multiplicity}, J. Algebraic Geom. 16 (2007), no. 1, 19--38.

\a \a  \a \small {{\it E.Ballico, Dept. of Mathematics, University of Trento, Povo (TN), Italy, email: ballico@science.unitn.it}} 
\b {\it M.Id\`a, Dip. di Matematica, Universit\`a di Bologna, Italia, e-mail: ida@dm.unibo.it}

\end{document}